\documentclass[12pt,a4paper,draft]{article}
\usepackage[english]{babel}
\usepackage{amssymb, amsmath, latexsym, amsfonts, amsthm}
\begin{document}
\renewcommand{\refname}{References}
\newtheorem{theorem}{Theorem}
\newtheorem{lemma}{Lemma}
\newtheorem{corollary}{Corollary}
\begin{center}
{\bf On a Fr\'echet space of entire functions rapidly decreasing on the real line} \footnote {
{\it  Key words}: Fourier transform, entire functions, convex functions. 

{\it AMS  Subject Class}. (2010): 30D15, 42A38}
\end{center}
\begin{center}
M. I. MUSIN
\end{center}
\renewcommand{\abstractname}{}
\begin{abstract}
{\sc Abstract}.
A weighted space of entire functions rapidly decreasing on the real line is considered in the paper. A 
growth of these functions along the imaginary axis is controlled by some system of weight functions. 
The Fourier transform of functions of this space is studied. 
Equivalent description of the considered space in terms of estimates 
on derivatives of functions on real line is obtained.   
\end{abstract}

\section {Introduction} 

\hspace{20pt}

For $u \in {\mathbb R}^n \ ({\mathbb C}^n)$ let $\Vert u \Vert$ be the Euclidean norm of $u$ in 
${\mathbb R}^n \ ({\mathbb C}^n)$. 

Let $\varPhi$ be a family of continuous functions on ${\mathbb R}^n$ such that:

${\varPhi 1}$. for all $\varphi_1, \varphi_2 \in \varPhi$ there exists $\varphi_3 \in \varPhi$ such that 
$$
\min (\varphi_1(x), \varphi_2(x)) - \varphi_3(x) \to +\infty, \ x \to \infty;
$$

${\varPhi 2}$. for each $\varphi \in \varPhi$
$$
\displaystyle \lim_{x \to  \infty} \frac {\varphi(x)}{\Vert x \Vert}= + \infty.
$$

For arbitrary $\varphi \in \varPhi$ and $k \in {\mathbb Z}_+$ let
$$
S_k(\varphi) = \{f \in H({\mathbb C}^n): 
q_{\varphi, k}(f) = \sup_{z \in {\mathbb C}^n} 
\frac 
{\vert f(z)\vert (1 + \Vert z \Vert)^k}
{e^{\varphi (Im z)}} < \infty \}.
$$
Here $H({\mathbb C}^n)$ is a space of entire functions on ${\mathbb C}^n$. 
Let 
$
S(\varPhi)= \bigcap \limits_{\varphi \in \varPhi, k \in {\mathbb Z_+}} S_k(\varphi).
$
With usual operations of addition and multiplication by complex numbers $S(\varPhi)$ is a linear space. Endow $S(\varPhi)$ with a topology defined by the system of norms $q_{\varphi, k}$ ($\varphi \in \varPhi$, $k \in {\mathbb Z}_+$). 

If $\varPhi$ consists of functions $\varphi (\varepsilon \Vert x \Vert) \ (\varepsilon > 0)$, 
where
$\varphi$ is a nonnegative continuous nondecreasing function  on $[0, \infty)$ such that: 

1). $\displaystyle \lim_{t \to + \infty} \frac {\varphi(t)}{t}= + \infty$;

2). the function $\psi(t) = \varphi(e^t)$ is convex on $[0, \infty)$;

3). there exist numbers $h > 1$  и $K > 0$ such that 
$$
2\varphi(t) \le \varphi(ht) + K, \ t \in [0, \infty),
$$
then the space $S(\varPhi)$ was considered for $n=1$ in \cite {Marat} and for $n \ge 1$ in \cite {MZ}. In \cite {Marat} and \cite {MZ} the space  $S(\varPhi)$ was described in terms of estimates on derivatives (partial derivatives) of functions on real line (on ${\mathbb R}^n$) and Fourier transformation of functions of $S(\varPhi)$ was studied.

In this paper we continue to study Fr\'echet spaces of entire functions rapidly decreasing on the real line and such that their growth along the imaginary axis is majorized with a help of some family of weight functions which are not necessary convex on the real line and satisfy more restrictive conditions than in \cite {Marat} and \cite {MZ}. Namely, 
let $\varphi$ be a nonnegative nondecreasing continuous function on $[0, \infty)$ such that: 

1). $\displaystyle \lim_{x \to + \infty} \frac {\varphi(x)}{x}= + \infty$;

2). $\forall h > 1 \ \exists K_h > 0$ \ 
$
2\varphi(x) \le \varphi(hx) + K_h, \ x \in [0, \infty);
$

3). $\varphi(e^x)$ is convex on $[0, \infty)$.

For example, the function $e^{x^{\alpha}}$ with $\alpha > 0$ satisfies these conditions. 

Let $\sigma > 0$. 
For arbitrary $\varepsilon > 0$ and $k \in {\mathbb Z}_+$ let
$$
S_{\varepsilon, k}(\varphi) = \{f \in H({\mathbb C}): 
p_{\varepsilon, k}(f) = \sup_{z \in {\mathbb C}} 
\frac 
{\vert f(z)\vert (1 + \vert z \vert)^k}
{e^{\varphi ((\sigma + \varepsilon) \vert Im \ z \vert)}} < \infty \}.
$$
Let 
$
S_{\sigma}(\varphi)= \bigcap \limits_{\varepsilon > 0, k \in {\mathbb Z_+}} S_{\varepsilon, k}(\varphi).
$
With usual operations of addition and multiplication by complex numbers $S_{\sigma}(\varphi)$ is a linear space. 
The family of norms $p_{\varepsilon, k}$ ($ \varepsilon> 0$, $k \in {\mathbb Z}_+$) defines a locally convex topology in $S_{\sigma}(\varphi)$. Endowed with this topology $S_{\sigma}(\varphi)$ is a Fr\'echet space. 

The aim of the article is to describe the space  $S_{\sigma}(\varphi)$ in terms of estimates 
on derivatives of functions on real line and to study Fourier transformation of functions of $S_{\sigma}(\varphi)$. 

{\bf 1.2. Main results}. For an arbitrary function $g$ on $[0, \infty)$ let $g[e]$ be a function on $[0, \infty)$ defined by the rule $g[e](x) = g(e^x), \ x \ge 0$. 

For a continuous function $g: [0, \infty) \to {\mathbb R}$  such that
$
\displaystyle \lim_{x \to + \infty} \frac {g(x)}{x}= + \infty,
$
let 
$
g^*(x) = \displaystyle \sup \limits_{y >0}(xy - g(y)), \ x \ge 0. 
$
The function $g^*$ is called Young conjugate of $g$ \cite {Ev}. 
Note that
$
\displaystyle \lim_{x \to + \infty} \frac {g^*(x)}{x}= + \infty.
$
It is well known that if $g$ is convex then $(g^*)^*=g$ \cite {Ev}.

Let $\psi$ be a function defined as follows: $\psi(x) = \varphi(e^x), \ x \in~{\mathbb R}$.

The following two theorems (proved in section 3) give another description of the space $S_{\sigma}(\varphi)$.

\begin{theorem}
Let $f \in S_{\sigma}(\varphi)$. 
Then 
$\forall \varepsilon> 0 \ \forall m \in {\mathbb Z}_+ \
\exists c_{\varepsilon, m}>0$ \ $\forall n~\in~{\mathbb Z}_+ \ \forall x \in {\mathbb R}$
$$
\vert x^m f^{(n)}(x)\vert \le c_{\varepsilon, m} {n}! 
(\sigma + \varepsilon)^n e^{-\psi^*(n)} \ .
$$
\end{theorem}

\begin{theorem}
Let $f \in C^{\infty}({\mathbb R})$ and for all $\varepsilon>0$ and $m \in {\mathbb Z}_+$ there exists 
a number $d_{\varepsilon, m} > 0$ such that for each $n \in {\mathbb Z}_+$
$$
(1+ \vert x \vert)^m \vert f^{(n)}(x) \vert \le  
d_{\varepsilon, m} (\sigma + \varepsilon)^n n! e^{-\psi^*(n)}, \ x \in {\mathbb R}.
$$
Then $f$ admits (the unique) extension to entire function belonging to $S_{\sigma}(\varphi)$.
\end{theorem}

For $\varepsilon>0, m \in {\mathbb Z_+}$ let 
$$
G_{\varepsilon, m}(\psi^*) = \{f \in C^m({\mathbb R}): 
\Vert f \Vert_{\varepsilon, m} = \max_{0 \le n \le m} 
\sup \limits_{x \in {\mathbb R}, k \in {\mathbb Z_+}}
\frac 
{\vert x^k f^{(n)}(x) \vert}
{k! (\sigma + \varepsilon)^k 
e^{-\psi^*(k)}} < \infty \}.
$$
Let
$
G_{\sigma}(\psi^*)= \bigcap \limits_{\varepsilon > 0, m \in {\mathbb Z_+}}G_{\varepsilon, m}(\psi^*).
$
With usual operations of addition and multiplication by complex numbers $G_{\sigma}(\psi^*)$ is a linear space. 
Endow $G_{\sigma}(\psi^*)$ with a topogy defined by the family of norms $\Vert f \Vert_{\varepsilon, m}$ ($\varepsilon >~0$, $m \in {\mathbb Z}_+$).

Fourier transform $\tilde f$ of $f \in S_{\sigma}(\varphi)$ is defined by the formula
$$
\tilde f(x) = \int_{{\mathbb R}} f(\xi) e^{-i x \xi} \ d \xi , \ x \in {\mathbb R}.
$$

In section 4 the following theorem is proved.
 
\begin{theorem}
Fourier transform establishes an isomorphism of spaces $S_{\sigma}(\varphi)$ and $G_{\sigma}(\psi^*)$.
\end{theorem}

If $\varphi$ is convex on $[0, \infty)$ then the space $G_{\sigma}(\psi^*)$ admits more simple description  (see section 5). Namely, the following theorem holds.

\begin{theorem}
Let $\varphi$ be convex on $[0, \infty)$. Then the space $G_{\sigma}(\psi^*)$ consists of functions $f \in C^{\infty}({\mathbb R})$ such that for each $\varepsilon > 0$ and  $n \in {\mathbb Z}_+$ there exists a number  $C_{\varepsilon,n}>0$ such that
$$
\vert f^{(n)}(x) \vert \le C_{\varepsilon, n} e^{-\varphi^*(\frac {\vert x \vert}
{\sigma + \varepsilon})}, \ x \in {\mathbb R}.
$$
\end{theorem}

\section{Auxiliary results}

Note that the condition 2) on $\varphi$ is equivalent to the following condition on $\psi$: 
$
\forall h > 1 \ \exists K_h > 0
$
$$
2\psi(x) \le \psi(x+\ln h) + K_h,\ x \in {\mathbb R}.
$$

\begin{lemma} 
For each $M>0$ there exists a number $A_M>0$ such that
$$
\psi^*(x) \le x \ln\frac {x}{M} - x + A_M, \ x > 0.
$$ 
\end{lemma}

{\bf Proof}. By definition of $\psi$ and the condition 1) on $\varphi$ it follows that for each $M>0$ there exists a number $A_M>0$ such that for all $y \ge 0$ \ 
$
\psi(y) \ge M e^y - A_M.
$
Hence, 
$$
\psi^*(x)=\sup_{y >0}(xy -\psi(y)) \le \sup_{y >0}(xy -M e^y) + A_M \le 
$$
$$
\le \sup_{y \in {\mathbb R}}(xy -M e^y) + A_M = 
x \ln\frac {x}{M} - x + A_M.
$$ 

From Lemma 1 we get the following

\begin{corollary} 
For each $b>0$ the series \ $\displaystyle \sum_{j=0}^{\infty} \frac {e^{\psi^*(j)}}{b^j j!}$ \ converges.
\end{corollary}

\begin{lemma} 
Let $g$ be a convex continuous function on $[0, \infty)$ such that
$
\displaystyle \lim_{x \to + \infty} \frac {g(x)}{x}= + \infty.
$
Then for each $\varepsilon > 0$ there exists a number $B_{\varepsilon} > 0$ such that
\begin{equation}
2 g(x) \le g(x+\varepsilon) + B_{\varepsilon}, \ x \ge 0,
\end{equation}
iff for each $\varepsilon > 0$ there exists a number $C_{\varepsilon} > 0$ such that 
\begin{equation}
g^*(x+y) \le g^*(x) + g^*(y) + \varepsilon (x+y) + C_{\varepsilon}, \ x, y \ge 0.
\end{equation}
\end{lemma} 

{\bf Proof}. {\bf Necessity}. Note first that
\begin{equation}
g^*(x) \ge -\inf \limits_{\xi \ge 0} g(\xi), \ x \ge 0.
\end{equation}
Let $\varepsilon > 0$ be arbitrary. Obviously, for all $x, y, t \in [0, \infty)$ we have
$$
g^*(x) + g^*(y) \ge (x+y) t - 2 g(t).
$$
In view of (1) for all  $x, y, t \ge 0$
$$
g^*(x) + g^*(y) \ge (x+y) (t+\varepsilon) - g(t+\varepsilon) - B_{\varepsilon} - \varepsilon(x+y).
$$
Therefore, for all $x, y \in [0, \infty)$
\begin{equation}
g^*(x) + g^*(y) \ge \sup_{\xi \ge \varepsilon} ((x+y)\xi - g(\xi)) - B_{\varepsilon} - \varepsilon(x+y).
\end{equation}
Further, for all $x, y \in [0, \infty)$
$$
\displaystyle\sup_{0 \le t < \varepsilon} ((x+y) \xi - g(\xi)) \le (x+y) \varepsilon  - \inf \limits_{0 \le \xi < \varepsilon} g(\xi) \le 
(x+y) \varepsilon - \inf \limits_{\xi \ge 0} g(\xi).
$$
Taking into account (3) we have 
$$
\displaystyle\sup_{0 \le t < \varepsilon} ((x+y) \xi - g(\xi)) \le 
(x+y) \varepsilon + g^*(x) \le 
$$
$$
\le (x+y) \varepsilon + g^*(x) + g^*(y) + \inf \limits_{\xi \ge 0} g(\xi).
$$
From this and inequality (4) we obtain 
$$
g^*(x+y) \le g^*(x) + g^*(y) + \varepsilon (x+y) + C_{\varepsilon}, \ x, y \ge 0.
$$
where $C_{\varepsilon}=\max (B_{\varepsilon}, \inf \limits_{\xi \ge 0} g(\xi))$.

{\bf Sufficiency}. Recall that by the inversion formula for Young transform \cite {Ev} $g=(g^*)^*$. 
Using this and (2) we have for each $\varepsilon > 0$
$$
2 g(x) =  \displaystyle \sup_{u \ge 0} (2 xu - 2 g^*(u)) 
\le 
\displaystyle \sup_{u \ge 0} (2 xu - g^*(2u) + 2 \varepsilon u + C_{\varepsilon}) = 
$$
$$ 
= \displaystyle \sup_{u \ge 0} ((x+\varepsilon) t - g^*(t)) + C_{\varepsilon}= g(x+\varepsilon) + C_{\varepsilon}.
$$
Put $B_{\varepsilon}=C_{\varepsilon}$. Thus, we got the inequality (1). The proof is complete.

The space $G_{\sigma}(\psi^*)$ admits another description. 
For $\varepsilon>0$, $m \in {\mathbb Z}_+$ let 
$$
Q_{\varepsilon, m}(\psi^*) = \{f \in C^m({\mathbb R}): 
s_{\varepsilon, m}(f) = \max_{0 \le n \le m}
\sup_{x \in {\mathbb R}, k  \in {\mathbb Z}_+}  
\frac 
{(1 + \vert x \vert)^k \vert f^{(n)}(x) \vert}
{k! (\sigma + \varepsilon)^k 
e^{-\psi^*(k)}} < \infty \}.
$$
Let $Q_{\sigma}(\psi^*)$ be a projective limit of the spaces $Q_{\varepsilon, m}(\psi^*)$.

The following lemma holds.

\begin{lemma} 
$Q_{\sigma}(\psi^*)=G_{\sigma}(\psi^*)$. 
\end{lemma} 

{\bf Proof}. Let $f \in Q_{\sigma}(\psi^*)$. 
ТThen for all $\varepsilon > 0, m \in {\mathbb Z}_+$  we have 
\begin{equation}
\Vert f \Vert_{\varepsilon, m} 
\le s_{\varepsilon, m}(f).
\end{equation}
This means that $f \in G_{\sigma}(\psi^*)$. Moreover, the embedding mapping $I: Q_{\sigma}(\psi^*) \to G_{\sigma}(\psi^*)$ is continuous. 

Now let $f \in G_{\sigma}(\psi^*)$, $\varepsilon > 0$ and $m \in {\mathbb Z}_+$ are arbitrary. 
Then 
$\Vert f \Vert_{\frac \varepsilon 2, m} < \infty$. 
Hence, for $n \in {\mathbb Z_+}$ such that $0 \le n \le m$ we have
\begin{equation}
\vert f^{(n)}(x) \vert \le \Vert f \Vert_{\frac \varepsilon 2, m} e^{-\psi^*(0)}, \ x \in {\mathbb R}.
\end{equation}
Choose $\delta=\delta(\varepsilon) > 0$ such small that
$(1+\delta)(\sigma + \frac {\varepsilon}{2}) \le \sigma + \varepsilon$. 
Note that for $n \in {\mathbb Z_+}$ such that $0 \le n \le m$ 
$$
\sup_{\vert x \vert \le \frac {1}{\delta}, k  \in {\mathbb Z}_+}  
\frac 
{(1 + \vert x \vert)^k \vert f^{(n)}(x) \vert}
{k! (\sigma + \varepsilon)^k 
e^{-\psi^*(k)}} \le 
\sup_{\vert x \vert \le \frac {1}{\delta}, k  \in {\mathbb Z}_+}  
\frac 
{(1+\frac {1}{\delta})^k \vert f^{(n)}(x) \vert}
{k! (\sigma + \varepsilon)^k 
e^{-\psi^*(k)}} \ .
$$
Since for each $b > 0$ \ 
$
\lim \limits_{k \to \infty} \frac {e^{\psi^*(k)}} {k! b^k }=0,
$
then there is a number $C(\varepsilon) > 1$ such that for $n \in {\mathbb Z_+}$ such that $0 \le n \le m$ we have
$$
\sup_{\vert x \vert \le \frac {1}{\delta}, k  \in {\mathbb Z}_+}  
\frac 
{(1 + \vert x \vert)^k \vert f^{(n)}(x) \vert}
{k! (\sigma + \varepsilon)^k 
e^{-\psi^*(k)}} \le C(\varepsilon) \sup_{\vert x \vert \le \frac {1}{\delta}}  \vert f^{(n)}(x) \vert  .
$$
Taking into account (6) we have  
\begin{equation}
\max_{0 \le n \le m} \sup_{\vert x \vert \le \frac {1}{\delta}, k  \in {\mathbb Z}_+}  
\frac 
{(1 + \vert x \vert)^k \vert f^{(n)}(x) \vert}
{k!(\sigma + \varepsilon)^k 
e^{-\psi^*(k)}} \le C(\varepsilon) \Vert f \Vert_{\frac \varepsilon 2, m}e^{-\psi^*(0)}.
\end{equation}
Since for $\vert x \vert > \frac {1}{\delta}$ \ 
$1 + \vert x \vert < (1+\delta) \vert x \vert$ then for $n \in {\mathbb Z_+}$ such that $0 \le n \le m$ 
we have
\begin{equation}
\sup_{\vert x \vert > \frac {1}{\delta}, k  \in {\mathbb Z}_+}  
\frac 
{(1 + \vert x \vert)^k \vert f^{(n)}(x) \vert}
{k! (\sigma + \varepsilon)^k e^{-\psi^*(k)}} \le 
\sup_{\vert x \vert > \frac {1}{\delta}, k  \in {\mathbb Z}_+}  
\frac 
{((1+\delta)\vert x \vert)^k \vert f^{(n)}(x) \vert}
{k! (\sigma + \varepsilon)^k e^{-\psi^*(k)}} \le \Vert f \Vert_{\frac \varepsilon 2, m}.
\end{equation}
From the estimates (7) and (8) it follows that for each $\varepsilon>0$ we have
\begin{equation}
s_{\varepsilon, m}(f) \le C_1(\varepsilon) \Vert f \Vert_{\frac \varepsilon 2, m}, \ f \in G_{\sigma}(\psi^*),
\end{equation}
where $C_1(\varepsilon)=\max (1, C(\varepsilon) e^{-\psi^*(0)})$.

From the estimates (5) and (9) it follows that $Q_{\sigma}(\psi^*)=G_{\sigma}(\psi^*)$ topologically.

\section{Equivalent description of the space $S_{\sigma}(\varphi)$}

{\bf Proof of Theorem 1}. Let $f \in S_{\sigma}(\varphi)$. Using Cauchy integral formula we have for all 
$m, n \in {\mathbb Z}_+$ and $R>0$
$$
(1 + \vert x \vert)^m f^{(n)}(x) = 
\frac {n! }{2\pi i} 
\displaystyle \int_{L_R(x)}
\frac 
{(1 + \vert x \vert)^m f(\zeta)}
{(\zeta - x)^{n +1}} \ d \zeta , \ x \in {\mathbb R},
$$
where
$L_R(x)= \{\zeta \in {\mathbb C}: \vert \zeta - x \vert = R \}$.
From this for each $\varepsilon>0$ we have
$$
 (1 + \vert x \vert)^m \vert f^{(n)}(x) \vert 
\le 
n! 
\max_{\zeta \in L_R} \frac {(1 + \vert  \zeta - x \vert)^m (1 + \vert  \zeta \vert)^m \vert f(\zeta) \vert} {R^{n}} 
\le 
$$
$$
\le 
n! p_{\frac {\varepsilon}{2}, m}(f)  
\frac 
{(1 + R)^m e^{\varphi ((\sigma + \frac {\varepsilon}{2}) R)}} 
{R^{n}} \ .
$$
Using the condition 2) on $\varphi$ and the fact that $\varphi$ is nondecreasing on $[0, \infty)$ we have for some  $d_{\varepsilon, m}>~0$ for each $R>0$
$$
(1 + \vert x \vert)^m \vert f^{(n)}(x) \vert 
\le 
d_{\varepsilon, m} n!p_{\frac {\varepsilon}{2}, m}(f) 
\frac 
{e^{\varphi ((\sigma + \varepsilon)R)}} {R^n} \ .
$$
From this for all $x \in {\mathbb R}$ we have
$$
(1 + \vert x \vert)^m \vert f^{(n)}(x) \vert 
\le 
d_{\varepsilon, m} n! p_{\frac {\varepsilon}{2}, m}(f) 
(\sigma + \varepsilon)^n \exp({-\sup_{R \ge 1} (n \ln R - \varphi (R))}) = 
$$
$$
=d_{\varepsilon, m} n! p_{\frac {\varepsilon}{2}, m}(f)  
(\sigma + \varepsilon)^n
\exp({-\sup_{r \ge 0}
(n r  - \psi (r))}) =
$$
$$
=
d_{\varepsilon, m} n! p_{\frac {\varepsilon}{2}, m}(f)  
(\sigma + \varepsilon)^n e^{-\psi^*(n)}.
$$

Theorem 1 is proved.

{\bf Proof of Theorem 2}. Let $f \in C^{\infty}({\mathbb R})$ and numbers 
$\varepsilon>0, m \in {\mathbb Z}_+$ are arbitrary. By the assumption there exist numbers $d > 0$ (depending on  $\varepsilon>0$) and $m \in {\mathbb Z}_+$ such that for each $n \in {\mathbb Z}_+$
\begin{equation}
(1+ \vert x \vert)^m \vert f^{(n)}(x) \vert \le  
d (\sigma + \frac {\varepsilon}{2})^n n! e^{-\psi^*(n)}, \ x \in {\mathbb R}.
\end{equation}
In particular, 
$
\vert f^{(n)}(x) \vert \le d (\sigma + \frac {\varepsilon}{2})^n n!e^{-\psi^*(n)}, \ x \in {\mathbb R}.
$
Taking into account that
$
\displaystyle \lim_{x \to + \infty} \frac {\psi^*(x)}{x}= + \infty,
$
it is easy to see that the sequence
$
\left(\sum_{n=0}^k 
\frac 
{f^{(n)}(0)}{n!} x^n\right)_{k=0}^{\infty}
$
converges to $f$ informly on compacts of the real line and the series
$
\displaystyle 
\sum_{n=0}^{\infty} 
\frac 
{f^{(n)}(0)}{n!} z^n
$
converges informly on compacts of ${\mathbb C}$ and, hence, its sum $F_f(z)$ is an entire function in ${\mathbb C}$. 
Note that $F_f(x) = f(x), \ x \in~{\mathbb R}$. The uniqueness of holomorphic continuation is obvious.

Using the equality
$
F_f(z) = \displaystyle 
\sum_{n=0}^{\infty} 
\frac 
{f^{(n)}(x)}{n!} (iy)^n, \ z = x+iy \ (x,y \in {\mathbb R}), 
$
and the inequality (10) we can estimate a growth  of $F_f$. Choose $\delta=\delta(\varepsilon) > 0$ so small that 
$(1+\delta)(\sigma + \frac {\varepsilon}{2}) \le \sigma + \varepsilon$. 
For all $m \in {\mathbb Z}_+$ and $\varepsilon > 0$  we have
$$
(1 + \vert z \vert)^m \vert F_f(z) \vert \le \sum_{n=0}^{\infty}  \frac 
{(1 + \vert x \vert)^m (1 + \vert y \vert)^{m+n} \vert f^{(n)}(x) \vert }{n!} 
\le 
$$
$$
\le 
\sum_{n=0}^{\infty}
\frac {d (\sigma + \frac {\varepsilon}{2})^n}
{e^{\psi^*(n)}} (1 + \vert y \vert)^{n + m} \le 
$$
$$
\le
d(1 + \vert y \vert)^m
\sum_{n=0}^{\infty}
\left(\frac {1}{1 + \delta}\right)^n
e^{\sup \limits_{n \ge 0} (n \ln ((1 + \delta) (\sigma + \frac {\varepsilon}{2}) (1 + \vert y \vert)) - \psi^*(n))} \le
$$
$$
\le \frac {d(1 + \delta)} {\delta} 
(1 + \vert y \vert)^m
e^{\sup \limits_{\xi \ge 0} (\xi \ln ((\sigma + \varepsilon) (1 + \vert y \vert)) - \psi^*(\xi))} = 
$$
$$
= 
\frac {d(1 + \delta)} {\delta} 
(1 + \vert y \vert)^m e^{\psi(\ln ((\sigma + \varepsilon) (1 + \vert y \vert)))}.
$$
Thus, 
\begin{equation}
(1 + \vert z \vert)^m \vert F_f(z) \vert \le \gamma_{\varepsilon} d
(1 + \vert y \vert)^m 
e^{\varphi ((\sigma + \varepsilon) (1 + \vert y \vert))}, \ z \in {\mathbb C},
\end{equation}
where $\gamma_{\varepsilon} > 0$ is some constant (depending on $\varepsilon$). 
Using the condition 2) on nondecreasing  function $\varphi$ it is possible to find a number $C_{\varepsilon, m, \varphi, \sigma} >0$ such that
$$
(1 + \vert z \vert)^m \vert F(z) \vert \le C_{\varepsilon, m, \varphi, \sigma} d e^{\varphi ((\sigma + 2 \varepsilon) \vert y \vert)}, \ z \in  {\mathbb C}. 
$$
Therefore,  $F_f \in S_{\sigma}(\varphi)$. 

Theorem 2 is proved.

\section{On Fourier transformation of functions of the space $S_{\sigma}(\varphi)$}

{\bf Proof of Theorem 3}. Let $f \in S_{\sigma}(\varphi)$. 
Then for all $\varepsilon > 0$ and for all $x \in {\mathbb R}$ we have
\begin{equation}
\vert \tilde f^{(n)}(x) \vert \le 
\int_{\mathbb R} 
\vert f(\xi) \vert
\vert \xi \vert^n \ d \xi 
\le 
\int_{\mathbb R} 
\frac 
{\vert f(\xi) \vert
(1 + \vert \xi \vert)^{n + 2}}{1 + \xi^2}
d \xi \le \pi p_{\varepsilon, n + 2}(f) \ .
\end{equation}
Since for all $m \in {\mathbb N}$, $n \in {\mathbb Z}_+$, $x, \eta \in {\mathbb R}$
$$
x^m \tilde f^{(n)}(x) = x^m \int_{\mathbb R} f(\zeta) 
(-i \zeta)^n 
e^{-i x \zeta} \ d \xi, \  \zeta = \xi + i\eta, 
$$
then
$$
\vert x^m \tilde f^{(n)}(x) \vert \le 
\int_{\mathbb R} 
\vert f(\zeta) \vert
\vert \zeta \vert^n 
e^{x \eta} \vert x \vert^m \ d \xi \le 
\int_{\mathbb R} 
\vert f(\zeta) \vert
(1 + \vert \zeta \vert)^{n + 2} 
e^{x \eta} \vert x \vert^m \ 
\frac {d \xi}{1 + \xi^2} \ .
$$

Consider the case $x \ne 0$. Put
$\eta = -\frac {x}{\vert x \vert} t, \ t > 0$. For all $t > 0$, $\varepsilon > 0$
$$
\vert x^m \tilde f^{(n)}(x) \vert \le \pi p_{\varepsilon, n + 2}(f) 
e^{-t \vert x \vert} 
e^{\varphi((\sigma + \varepsilon) t)} \vert x \vert^m  \le 
$$
$$
\le \pi p_{\varepsilon, n + 2}(f) 
e^{\sup \limits_{r>0} (-tr + m \ln r)} e^{\varphi((\sigma + \varepsilon) t)} \le
\pi p_{\varepsilon, n + 2}(f) 
e^{m \ln m - m - m \ln t} e^{\varphi((\sigma + \varepsilon) t)}.
$$
From this taking into account that
$$
\inf_{t>0} (- m \ln t + \varphi ((\sigma + \varepsilon)t)) = 
m \ln (\sigma + \varepsilon) - \sup_{u>0} (m \ln u - \varphi (u)) \le
$$
$$
\le
m \ln (\sigma + \varepsilon) - \sup_{u \ge 1} (m \ln u - \varphi (u))= m \ln (\sigma + \varepsilon) - \psi^*(m),
$$
we have
\begin{equation}
\vert x^m \tilde f^{(n)}(x) \vert \le \pi p_{\varepsilon, n + 2}(f) (\sigma + \varepsilon)^m  e^{m \ln m - m}
e^{- \psi^*(m)}.
\end{equation}
If $x=0$, then for $m \in {\mathbb N}$ and  
$n \in {\mathbb Z}_+$\ $x^m \tilde f^{(n)}(x) =0$.  From this, estimates (12) and (13) and since
$m^m \le e^m m!$ for all $m \in {\mathbb N}$, 
we have
for all $\varepsilon > 0, k \in {\mathbb Z}_+$ \
$$
\Vert \tilde f \Vert_{\varepsilon, k} \le \pi p_{\varepsilon, k + 2}(f), \ f \in S_{\sigma}(\varphi).
$$
This means that a linear 
transformation ${\mathcal F}: S_{\sigma}(\varphi) \to G_{\sigma}(\psi^*)$ acting by the rule $f \in S_{\sigma}(\varphi) \to \tilde f$ is continuous. 

Let us show that ${\mathcal F}$ is surjective. Let $g \in G_{\sigma}(\psi^*)$. Then (using Lemma 3) for all
$\varepsilon > 0$, $k, m \in {\mathbb Z}_+$ and $ n \in {\mathbb Z_+}$ such that
$n \le m$ we have
$$
{(1 + \vert x \vert)^k \vert g^{(n)}(x) \vert} \le 
s_{\varepsilon, m}(g)
(\sigma + \varepsilon)^k k! e^{-\psi^*(k)}, \ x \in {\mathbb R}.
$$
Let
$$
f(\xi) = 
\frac {1}{2 \pi} 
\int_{{\mathbb R}} g(x) e^{i x \xi} \ dx, \ 
\xi \in {\mathbb R}.
$$
For each $n \in {\mathbb Z_+}$
$$
f^{(n)}(\xi) = 
\frac {1}{2 \pi} 
\int_{{\mathbb R}} g(x) (ix)^n e^{i x \xi} \ dx, \ 
\xi \in {\mathbb R}.
$$
From this (integrating by parts) for each $m \in {\mathbb Z_+}$
$$
(i\xi)^m f^{(n)}(\xi) = \frac {1}{2 \pi} (-1)^m 
\int_{{\mathbb R}} (g(x) (ix)^n)^{(m)} e^{i x \xi} \ dx, \ 
\xi \in {\mathbb R}.
$$
Let $r =\min (m, n)$. Then
$$
(i\xi)^m f^{(n)}(\xi) = \frac {1}{2 \pi} (-1)^m \int_{{\mathbb R}}
\displaystyle \sum_{j=0}^r C_m^j g^{(m-j)}(x) ((ix)^n)^{(j)} e^{ix \xi} \ dx , \ 
\xi \in {\mathbb R}.
$$
From this we have
$$
\vert \xi^m f^{(n)}(\xi)\vert \le 
\frac {1}{2\pi}
\displaystyle \sum_{j=0}^r
C_m^j \int_{{\mathbb R}} \vert g^{(m-j)}(x)
\vert \frac {n!}{(n-j)!} \vert x \vert^{n-j}  \ dx \le
$$ 
$$
\le
\frac {1}{2\pi}
\displaystyle \sum_{j=0}^r
C_m^j \frac {n!}{(n-j)!}
\int_{\mathbb R}
\vert g^{(m-j)}(x) \vert (1 + \vert x \vert)^{n  - j + 2} \ \frac {dx}{1 + x^2} \le
$$
$$
\le 
\frac {1}{2}
\displaystyle \sum_{j=0}^r
C_m^j \frac {n!}{(n-j)!}
s_{\varepsilon, m}(g) (n-j+2)! (\sigma + \varepsilon)^{n-j+2} e^{-\psi^*(n-j+2)} \le
$$
$$
\le
\frac {1}{2} n! (\sigma + \varepsilon)^{n+2} s_{\varepsilon, m}(g)
\displaystyle \sum_{j=0}^r
C_m^j (n-j+1)(n-j+2) (\sigma + \varepsilon)^{-j} e^{-\psi^*(n-j)} .
$$
Using conditions on $\psi$ and Lemma for each $\delta > 0$ it is possible to find a number $K_{\delta}>0$ such that
\begin{equation}
\psi^*(x+y) \le \psi^*(x) + \psi^*(y) + \delta (x+y) + K_{\delta}, \ x, y \ge 0.
\end{equation}
So for all $\xi \in {\mathbb R}$ and all $\varepsilon > 0$ and $\delta > 0$
$$
\vert \xi^m f^{(n)}(\xi)\vert \le 
\frac {1}{2} (n+2)! (\sigma + \varepsilon)^2 
((\sigma + \varepsilon) e^{\delta})^n 
s_{\varepsilon, m}(g) e^{-\psi^*(n)}
\displaystyle \sum_{j=0}^r 
\frac {e^{\psi^*(j) + K_{\delta}}}{(\sigma + \varepsilon)^j j!} \ .
$$
Choose $\delta=\delta(\varepsilon)$ so small that $(\sigma + \varepsilon) e^{\delta} < \sigma + 2 \varepsilon$ and  
put 
$$
c_{\varepsilon, m} = (\sigma + \varepsilon)^2 e^{K_{\delta}} \displaystyle \sum_{j=0}^{\infty} 
\frac {e^{\psi^*(j)}}{(\sigma + \varepsilon)^j j!} \ . 
$$
Then for all $n \in {\mathbb Z_+}$ and $\varepsilon > 0$ 
$$
\vert \xi^m f^{(n)}(\xi)\vert \le c_{\varepsilon, m} (n+1)(n+2)((\sigma + 2 \varepsilon))^n n! 
s_{\varepsilon, m}(g) e^{-\psi^*(n)}, \ \xi \in {\mathbb R}.
$$
Find a number $D_{\varepsilon} > 0$ such that for all $n \in {\mathbb Z_+}$
$$
(n+1)(n+2)(\sigma + 2\varepsilon)^n \le D_{\varepsilon}(\sigma + 3\varepsilon)^n.
$$
Then for all $n \in {\mathbb Z_+}$ and $\xi \in {\mathbb R}$
$$
(1 + \vert \xi \vert)^m \vert f^{(n)}(\xi)\vert \le 2^m D_{\varepsilon}(c_{\varepsilon, 0} + c_{\varepsilon, m}) 
s_{\varepsilon, m}(g) (\sigma + 3\varepsilon)^n n! e^{-\psi^*(n)}. 
$$
By Theorem 2 $f$ can be holomorphically continued (uniquely) to entire function belonging to $S_{\sigma}(\varphi)$. Denote this continuation by $F$. 
Obviously, $g = {\mathcal F}(F)$. 
In view of the inequalities (10) and (11) we have for some $\gamma_{\varepsilon} > 0$ 
$$
(1 + \vert z \vert)^m \vert F(z) \vert \le  \gamma_{\varepsilon} 2^m(c_{\varepsilon, 0} + c_{\varepsilon, m}) 
s_{\varepsilon, m}(g)(1 + \vert y \vert)^m e^{\varphi ((\sigma + 6 \varepsilon)\vert y \vert))}, \ z \in {\mathbb C}.
$$
Taking into account conditions 1) and 2) on $\varphi$ one can find a number $K_{\varepsilon, m, \varphi}~>~0$ 
such that\ 
$
p_{7\varepsilon, m}(f) \le K_{\varepsilon, m, \varphi} s_{\varepsilon, m}(g). 
$
From this and Lemma 3 we get that the inverse mapping ${\mathcal F}^{-1}$ is continuous. 

Thus, we have proved that Fourier transformation establishes a topological isomorphism between spaces  $S_{\sigma}(\varphi)$ and $G_{\sigma}(\psi^*)$.

\section{Special case of a function $\varphi$}

In the proof of Theorem 4 we will use the following result.

\begin{lemma} 
Let $g$ be a real-valued continuous function on  
$[0, \infty)$ such that
$
\displaystyle \lim_{x \to + \infty} \frac {g(x)}{x}= + \infty.
$
Then for each $\delta > 0$
$$
\displaystyle \lim_{x \to + \infty} \frac {g^*((1+\delta) x) - g^*(x)}{x}= + \infty.
$$
\end{lemma} 

{\bf Proof}. Let $\delta > 0$ be arbitrary. For each $x>0$ denote by $\xi(x)$ a point where the supremum  of the function $u_x(\xi)= x \xi - g(\xi)$ over $[0, \infty)$ is attained. 
Note that  $\xi(x) \to +\infty$ as $x \to +\infty$. Otherwise there are a number $M>0$ and a sequence $(x_j)_{j=1}^{\infty}$ of positive numbers $x_j$ converging to $+\infty$ such that $\xi(x_j) \le M$. 
Then $g^*(x_j)= x_j \xi(x_j) - g(\xi(x_j))$.  But it contradicts to the fact that 
$
\displaystyle \lim_{x \to + \infty} \frac {g^*(x)}{x}= + \infty.
$
Thus, 
$
\displaystyle \lim_{x \to + \infty} \xi(x)= + \infty.
$
From this and the inequality
$$
g^*((1+\delta) x)- g^*(x) \ge (1+\delta) x \xi(x) - g(\xi(x)) - x \xi(x) + g(\xi(x))=\delta x \xi(x), \ x>0,
$$
the assertion of lemma follows.

{\bf Proof of Theorem 4}. Let  $\varphi$ satisfies the conditions 1) -- 3) and be convex on $[0, \infty)$. 
Let us show that in this case  the space $G_{\sigma}(\psi^*)$ consists of functions $f \in C^{\infty}({\mathbb R})$ 
such that for all $\varepsilon > 0$, $n \in {\mathbb Z}_+$ there exists a number $C_{\varepsilon, n}>0$ such that
\begin{equation}
\vert f^{(n)}(x) \vert \le C_{\varepsilon, n} e^{-\varphi^*(\frac {\vert x \vert}
{\sigma + \varepsilon})}, \ x \in {\mathbb R}.
\end{equation}

Let us construct a continuous increasing convex function $\varphi_1$ on $[0, \infty)$ such that:

1. $\varphi_1(x) > 0$ for $x>0$;

2.  $\displaystyle \lim_{x \to 0+} \frac {\varphi_1(x)}{x}= 0$;

3. $\varphi_1$ coinsides with $\varphi$ out of the segment $[0, d]$ where $d > 0$ is some number.

Define a function $\psi_1$ on $[0, \infty)$ by the equality $\psi_1=\varphi_1[e]$. 

Obviously, there are positive numbers $s=s(\varphi, d)$ и $s_1 = s_1(\varphi, d)$ such that for all $x \ge 0$
\begin{equation}
\vert \varphi^*(x) - \varphi_1^*(x) \vert \le s, 
\end{equation}  
\begin{equation}
\vert \psi^*(x) - \psi_1^*(x) \vert \le s_1.
\end{equation}

Let $f \in G_{\sigma}(\psi^*)$ and $\varepsilon \in (0, 1)$ are arbitrary. 
By Lemma 3 for all $n, k \in {\mathbb Z}_+$ 
\begin{equation}
\vert f^{(n)}(x) \vert \le s_{\varepsilon, n}(f) 
\frac 
{k! (\sigma + \varepsilon)^k 
e^{-\psi^*(k)}}{(1 + \vert x \vert)^k}  \ , \ x \ge 0.
\end{equation}
Since for each $k \in {\mathbb N}$ \
$
k! < 3 \frac {k^{k+1}}{e^k},
$
then
\begin{equation}
\frac 
{k! (\sigma + \varepsilon)^k  e^{-\psi^*(k)}}{(1 + \vert x \vert)^k} < \frac {3(\sigma + \varepsilon)^k  k^{k+1} e^{-\psi^*(k)}} {e^k(1 + \vert x \vert)^k} \ , \ k \in {\mathbb N}.
\end{equation}
If $\sigma + \varepsilon \ge 1$, 
then from (19) 
we have for 
$k \in {\mathbb N}$, $t \in [k, k+1)$, $x \in {\mathbb R}$ and for all $\delta > 0$
$$
\frac 
{k! (\sigma + \varepsilon)^k  e^{-\psi^*(k)}}{(1 + \vert x \vert)^k} 
\le
\frac 
{3 (\sigma + \varepsilon)^t t^{t+1} e^{-\psi^*(t) + \psi^*(1) + \delta (t+1) + K_{\delta} + 1}}{e^t(1 + \vert x \vert)^t} (1 + \vert x \vert) = 
$$
$$
=
3e^{K_{\delta} + 1 + \delta + \psi^*(1)} e^{t \ln (\sigma + \varepsilon) +(t+1) \ln t -\psi^*(t) + \delta t - t - t \ln (1 + \vert x \vert)}(1 + \vert x \vert).
$$
If $\sigma + \varepsilon < 1$, then from (19) (taking into account the inequality (14) and that $\psi^*$ is nondecreasing) we have for $k \in {\mathbb N}$, $t \in [k, k+1)$, $x \in {\mathbb R}$ and for all $\delta > 0$
$$
\frac 
{k! (\sigma + \varepsilon)^k  e^{-\psi^*(k)}}{(1 + \vert x \vert)^k} \le 
\frac 
{3 (\sigma + \varepsilon)^t t^{t+1} e^{-\psi^*(t) + \psi^*(1) + \delta (t+1) + K_{\delta} + 1}}{(\sigma + \varepsilon)e^t(1 + \vert x \vert)^t} (1 + \vert x \vert) = 
$$
$$
=\frac {3e^{K_{\delta} + 1 + \delta + \psi^*(1)}}{\sigma + \varepsilon} e^{t \ln (\sigma + \varepsilon) +(t+1) \ln t -\psi^*(t) + \delta t - t - t \ln (1 + \vert x \vert)}(1 + \vert x \vert).
$$
Put $C_{\sigma, \varepsilon, \delta}=3e^{K_{\delta} + 1 + \delta + \psi^*(1)}$ if $\sigma + \varepsilon \ge 1$, 
$C_{\sigma, \varepsilon, \delta}=\frac {3e^{K_{\delta} + 1 + \delta + \psi^*(1)}}{\sigma + \varepsilon}$ if $\sigma + \varepsilon < 1$.
Thus, for all $\varepsilon > 0$ and $\delta > 0$ we have
$$
\frac 
{k! (\sigma + \varepsilon)^k  e^{-\psi^*(k)}}{(1 + \vert x \vert)^k} \le C_{\sigma, \varepsilon, \delta}
e^{t \ln (\sigma + \varepsilon) +(t+1) \ln t -\psi^*(t) + \delta t - t - t \ln (1 + \vert x \vert)}(1 + \vert x \vert).
$$
From this using (17)  and putting $K_{\sigma, \varepsilon, \delta}=C_{\sigma, \varepsilon, \delta} e^{s_1}$ we have
\begin{equation}
\frac 
{k! (\sigma + \varepsilon)^k  e^{-\psi^*(k)}}{(1 + \vert x \vert)^k} \le  K_{\sigma, \varepsilon, \delta} 
e^{t \ln (\sigma + \varepsilon) +(t+1) \ln t -\psi_1^*(t) + \delta t - t - t \ln (1 + \vert x \vert)}(1 + \vert x \vert).
\end{equation}

Note that if $u$ is an increasing convex function on $[0, \infty)$ such that
$$
\lim \limits_{x \to +\infty}
\displaystyle \frac {u(x)} {x} =
\lim \limits_{x \to 0+}
\displaystyle \frac {x} {u(x)} = + \infty,
$$
then the following formula holds
\begin{equation}
(u(e))^*(x) + (u^*(e))^*(x) = x \ln x - x, \
x > 0.
\end{equation}
This formula is a particular case of a general formula of changes of variables in Young transformation of convex operators  \cite {Kut}. It is noted in \cite {N-P}.

Using (21) we get from (20)
$$
\frac 
{k! (\sigma + \varepsilon)^k
e^{-\psi^*(k)}}{(1 + \vert x \vert)^k} < K_{\sigma, \varepsilon, \delta} (1 + \vert x \vert)
e^{t \ln \frac {(\sigma + \varepsilon) e^{\delta}}{1 + \vert x \vert} + \ln t + (\varphi_1^*[e])^*(t)} 
$$
Choose $\delta = \delta (\varepsilon)$ so small that 
$(\sigma + \varepsilon) e^{\delta} \le \sigma + \frac 3 2 \varepsilon$. 
Then
$$
\frac 
{k! (\sigma + \varepsilon)^k
e^{-\psi^*(k)}}{(1 + \vert x \vert)^k} < 
K_{\sigma, \varepsilon, \delta(\varepsilon)} (1 + \vert x \vert)
e^{t \ln \frac {\sigma + \frac 3 2 \varepsilon}{1 + \vert x \vert} + \ln t + (\varphi_1^*[e])^*(t)}
$$
Choose a number $M_{\sigma, \varepsilon} > 0$ so that 
$$
\frac 
{k! (\sigma + \varepsilon)^k
e^{-\psi^*(k)}}{(1 + \vert x \vert)^k} < M_{\sigma, \varepsilon} (1 + \vert x \vert)
e^{t \ln \frac {\sigma + 2 \varepsilon}{1 + \vert x \vert} + (\varphi_1^*[e])^*(t)} 
$$
From this it follows that
\begin{equation}
\inf_{k \in {\mathbb N}} 
\frac 
{k! (\sigma + \varepsilon)^k
e^{-\psi^*(k)}}{(1 + \vert x \vert)^k} \le M_{\sigma, \varepsilon} (1 + \vert x \vert)
\inf_{t \ge 1}
e^{t \ln \frac {\sigma + 2 \varepsilon}{1 + \vert x \vert} + (\varphi_1^*[e])^*(t)}.
\end{equation}
Obviously, 
\begin{equation}
\inf_{t \ge 1}
(t \ln \frac {\sigma + 2 \varepsilon}{1 + \vert x \vert} + (\varphi_1^*[e])^*(t)) \le 
\ln \frac {\sigma + 2 \varepsilon}{1 + \vert x \vert} + (\varphi_1^*[e])^*(1) \ .
\end{equation}
If
$\frac {\sigma + 2 \varepsilon}{1 + \vert x \vert} > 1$, then 
$$
\inf_{0 < t \le 1}
(t \ln \frac {\sigma + 2 \varepsilon}{1 + \vert x \vert} + (\varphi_1^*[e])^*(t)) \ge 
$$
$$
\ge 
\inf_{0 < t \le 1}
(t \ln \frac {\sigma + 2 \varepsilon}{1 + \vert x \vert}) + 
\inf_{0 < t \le 1}(\varphi_1^*[e])^*(t)) = (\varphi_1^*[e])^*(0).
$$
Using the inequality (23) we have 
$$
\inf_{0 < t \le 1}
(t \ln \frac {\sigma + 2 \varepsilon}{1 + \vert x \vert} + (\varphi_1^*[e])^*(t)) - (\varphi_1^*[e])^*(0) \ge 0 \ge 
$$
$$
\ge  
\inf_{t \ge 1}
(t \ln \frac {\sigma + 2 \varepsilon}{1 + \vert x \vert} + (\varphi_1^*[e])^*(t)) - \ln \frac {\sigma + 2 \varepsilon}{1 + \vert x \vert} - (\varphi_1^*[e])^*(1).
$$
Thus,
$$
\inf_{t \ge 1}
(t \ln \frac {\sigma + 2 \varepsilon}{1 + \vert x \vert} + (\varphi_1^*[e])^*(t)) \le 
$$
$$
\le
\inf_{0 < t \le 1}
(t \ln \frac {\sigma + 2 \varepsilon}{1 + \vert x \vert} + (\varphi_1^*[e])^*(t)) + \ln \frac {\sigma + 2 \varepsilon}{1 + \vert x \vert} + (\varphi_1^*[e])^*(1) - (\varphi_1^*[e])^*(0).
$$
Now it is clear that
$$
\inf_{t \ge 1}
(t \ln \frac {\sigma + 2 \varepsilon}{1 + \vert x \vert} + (\varphi_1^*[e])^*(t)) \le 
$$
$$
\le
\inf_{t > 0}
(t \ln \frac {\sigma + 2 \varepsilon}{1 + \vert x \vert} + (\varphi_1^*[e])^*(t)) + \ln \frac {\sigma + 2 \varepsilon}{1 + \vert x \vert} + (\varphi_1^*[e])^*(1) - (\varphi_1^*[e])^*(0).
$$
If
$\frac {\sigma + 2 \varepsilon}{1 + \vert x \vert} \le 1$, then
$$
\inf_{0 < t \le 1}
(t \ln \frac {\sigma + 2 \varepsilon}{1 + \vert x \vert} + (\varphi_1^*[e])^*(t)) \ge 
$$
$$
\ge 
\inf_{0 < t \le 1}
(t \ln \frac {\sigma + 2 \varepsilon}{1 + \vert x \vert}) + 
\inf_{0 < t \le 1}(\varphi_1^*[e])^*(t)) = -\ln \frac {\sigma + 2 \varepsilon}{1 + \vert x \vert} + (\varphi_1^*[e])^*(0).
$$
So using the inequality (23) we get
$$
\inf_{0 < t \le 1}
(t \ln \frac {\sigma + 2 \varepsilon}{1 + \vert x \vert} + 
(\varphi_1^*[e])^*(t)) + \ln \frac {\sigma + 2 \varepsilon}{1 + \vert x \vert} - (\varphi_1^*[e])^*(0) \ge 0 \ge 
$$
$$
\ge
\inf_{t \ge 1}
(t \ln \frac {\sigma + 2 \varepsilon}{1 + \vert x \vert} + (\varphi_1^*[e])^*(t)) - \ln \frac {\sigma + 2 \varepsilon}{1 + \vert x \vert} - (\varphi_1^*[e])^*(1).
$$
From this we obtain that
$$
\inf_{t \ge 1}
(t \ln \frac {\sigma + 2 \varepsilon}{1 + \vert x \vert} + (\varphi_1^*[e])^*(t)) \le 
$$
$$
\le
\inf_{0 < t \le 1}
(t \ln \frac {\sigma + 2 \varepsilon}{1 + \vert x \vert} + 
(\varphi_1^*[e])^*(t)) + 2 \ln \frac {\sigma + 2 \varepsilon}{1 + \vert x \vert} + (\varphi_1^*[e])^*(1)- (\varphi_1^*[e])^*(0).
$$
Obviously, 
$$
\inf_{t \ge 1}
(t \ln \frac {\sigma + 2 \varepsilon}{1 + \vert x \vert} + (\varphi_1^*[e])^*(t)) \le 
$$
$$
\le
\inf_{t > 0}
(t \ln \frac {\sigma + 2 \varepsilon}{1 + \vert x \vert} + 
(\varphi_1^*[e])^*(t)) + (\varphi_1^*[e])^*(1)- (\varphi_1^*[e])^*(0).
$$
Thus, in both cases we have
$$
\inf_{t \ge 1}
(t \ln \frac {\sigma + 2 \varepsilon}{1 + \vert x \vert} + (\varphi_1^*[e])^*(t)) \le 
$$
$$
\le
\inf_{t > 0}
(t \ln \frac {\sigma + 2 \varepsilon}{1 + \vert x \vert} + 
(\varphi_1^*[e])^*(t)) + \ln^+ \frac {\sigma + 2 \varepsilon}{1 + \vert x \vert} + (\varphi_1^*[e])^*(1)- (\varphi_1^*[e])^*(0),
$$
where $\ln^+t=\ln t$, если $t \ge 1$ и $\ln^+t=0$, если $t < 1$.
Going back to (22) we have
$$
\inf_{k \in {\mathbb N}} 
\frac 
{k! (\sigma + \varepsilon)^k
e^{-\psi^*(k)}}{(1 + \vert x \vert)^k} \le m_{\sigma, \varepsilon} (1 + \vert x \vert) 
e^{\inf \limits_{t > 0}
(t \ln \frac {\sigma + 2 \varepsilon}{1 + \vert x \vert} + 
(\varphi_1^*[e])^*(t)) + \ln^+ \frac {\sigma + 2 \varepsilon}{1 + \vert x \vert}},
$$
where
$m_{\sigma, \varepsilon}=M_{\sigma, \varepsilon} e^{(\varphi_1^*[e])^*(1)- (\varphi_1^*[e])^*(0)}$.
From this putting 
$c_{\sigma, \varepsilon}=\max(m_{\sigma, \varepsilon}, m_{\sigma, \varepsilon} (\sigma + 2 \varepsilon))$, 
we get
$$
\inf_{k \in {\mathbb N}} 
\frac 
{k! (\sigma + \varepsilon)^k
e^{-\psi^*(k)}}{(1 + \vert x \vert)^k} \le c_{\sigma, \varepsilon} 
e^{- \sup \limits_{t > 0}(t \ln \frac {1 + \vert x \vert} {\sigma + 2 \varepsilon} - (\varphi_1^*[e])^*(t))} 
(1 + \vert x \vert) =
$$
$$ 
=  c_{\sigma, \varepsilon}
e^{-(\varphi_1^*([e])^{**}(\ln \frac {1 + \vert x \vert} {\sigma + 2 \varepsilon})}(1 + \vert x \vert) = c_{\sigma, \varepsilon} 
e^{-\varphi_1^*([e])(\ln \frac {1 + \vert x \vert} {\sigma + 2 \varepsilon})}(1 + \vert x \vert) = 
$$
$$
= c_{\sigma, \varepsilon} 
e^{-\varphi_1^*(\frac {1 + \vert x \vert} {\sigma + 2 \varepsilon})}(1 + \vert x \vert).
$$
Using Lemma 4 we can find a number $\mu_{\sigma, \varepsilon} > 0$ such that
$$
\inf_{k \in {\mathbb N}} 
\frac 
{k! (\sigma + \varepsilon)^k 
e^{-\psi^*(k)}}{(1 + \vert x \vert)^k} \le \mu_{\sigma, \varepsilon} 
e^{-\varphi_1^*(\frac {1 + \vert x \vert} {\sigma + 3 \varepsilon})}.
$$
From this, (16) and (18) we have for each positive number $\varepsilon$ that
$$
\vert f^{(n)}(x) \vert \le s_{\varepsilon, n}(f)
\mu_{\sigma, \varepsilon} e^s
e^{-\varphi^*(\frac {\vert x \vert} {\sigma + 3 \varepsilon})}, \ x \in {\mathbb R}.
$$
Thus, we got the estimate of the form (15).

Now let $f \in C^{\infty}({\mathbb R})$ satisfies the estimate of the form (15).
Using (16) for each $\varepsilon  > 0$ и $n \in {\mathbb Z_+}$ we can find a number 
$R_{\varepsilon, n} > 0$ such that
\begin{equation}
\vert f^{(n)}(x) \vert \le R_{\varepsilon, n}(f)
e^{-\varphi_1^*(\frac {\vert x \vert} {\sigma +  \varepsilon})}, \ x \in {\mathbb R}.
\end{equation}
Let us show that $f \in G_{\sigma}(\psi^*)$. 
For $x \ne 0$ and $n \in {\mathbb Z}_+$ the inequality (24) can be rewritten as follows:
$$
\vert f^{(n)}(x) \vert \le R_{\varepsilon, n}(f)  
e^{-\varphi_1^*[e]
(\ln  \frac 
{\vert x \vert} {\sigma + \varepsilon})}.
$$
Using the inversion formula for Young transformation we have
$$
\vert f^{(n)}(x) \vert \le R_{\varepsilon, n}(f)  
e^{- \sup \limits_{t > 0} (t \ln  \frac {\vert x \vert}{\sigma + \varepsilon} - (\varphi_1^*[e])^*(t))}, \ x \ne 0.
$$
Using the equality (21) we get
$$
\vert f^{(n)}(x) \vert \le R_{\varepsilon, n} 
e^{-\sup\limits_{t > 0} (t \ln  \frac {\vert e x \vert}{\sigma + \varepsilon} - t \ln t + \psi_1^*(t))}, \ x \ne 0.
$$
Thus, for all $\varepsilon  > 0, k \in  {\mathbb N}$
$$
\vert f^{(n)}(x) x^k \vert \le R_{\varepsilon, n} (\sigma + \varepsilon)^k \left(\frac {k}{e}\right)^k e^{- \psi_1^*(k)}, \ x \ne 0.
$$
Since $k^k \le e^k k!$ for $k \in {\mathbb N}$ then 
$$
\vert f^{(n)}(x) x^k \vert \le R_{\varepsilon, n} (\sigma + \varepsilon)^k k! e^{- \psi_1^*(k)}, \ k \in {\mathbb N}, x \ne 0.
$$
Obviously, this inequality holds at the point $x = 0$ for all $k \in  {\mathbb N}$  
and for all $x \in {\mathbb R}$ if $k=0$ (in view of (15)). Ising (17) we have for all $\varepsilon$ and $k \in {\mathbb N}$
$$
\vert f^{(n)}(x) x^k \vert \le R_{\varepsilon, n} e^s (\sigma + \varepsilon)^k k! e^{- \psi^*(k)}, \ k \in {\mathbb N}, \ x \in {\mathbb R}.
$$
Thus, $f \in G_{\sigma}(\psi^*)$ and Theorem 4 is proved.

M. I. Musin, 

Bashkirian state University, Department of Mathematics and Computer 

technologies, Z. Validy str., 32, Ufa, 450000, Russia

E-mail address: marat402@gmail.com


\begin{thebibliography}{99}

\bibitem {Marat} 
M. I. Musin, {\it On a space of entire functions fast decreasing on the real line}, Ufa Mathematical Journal, {\bf 4}:1 (2012), 129–-137. 

\bibitem {MZ} 
M. I. Musin, {\it On spaces of entire functions on ${\mathbb C}^n$ rapidly decreasing on ${\mathbb R}^n$}. Submitted to Matematicheskie Zametki (Mathematical Notes).

\bibitem {Ev}
M. A. Evgrafov, Asymptotic estimates and entire functions. Moscow, Nauka. 1979. 

\bibitem {Kut} 
S. S. Kutateladze, {\it Changes of variables in the Young transformation},  Soviet Mathematics. Doklady, {\bf  233}:6 (1977), 1039--1041. (in Russian) 

\bibitem{N-P}
V. V. Napalkov, S. V. Popyonov, {\it On Laplace transformation on weighted Bergman space of entire functions on ${\mathbb C}^n$}, Doklady Mathematics (Doklady Akademii Nauk), {\bf 55}:1, 110--112.

\end{thebibliography}
\end{document}